\documentclass[12pt]{article}
\usepackage{amsfonts}
\usepackage{amssymb}
\usepackage{amsthm}
\usepackage{newlfont}

\makeatletter
\newfont{\footsc}{cmcsc10 at 8truept}
\newfont{\footbf}{cmbx10 at 8truept}
\newfont{\footrm}{cmr10 at 10truept}

 \newtheorem{thm}{Theorem}
 \newtheorem{cor}[thm]{Corollary}
 
 \newtheorem{dfn}[thm]{Definition}
 \newtheorem{con}[thm]{Conjecture}
 \newtheorem{opn}[thm]{Open Question}

\title{The Pi-Pebbling Function}

\author{T. Baillie Arnold\thanks{Partially Supported by NSF Grant DMS-0243774}\\
\small Department of Mathematics\\[-0.8ex]
\small Bowdoin College, Brunswick, Maine\\[-0.8ex]
\small \texttt{tarnold2@bowdoin.edu}}

\date{\small
Submitted: June 26, 2005;\\  
\small Mathematics Subject Classifications: 05C99, 05C35}

\begin{document}
\maketitle

\begin{abstract}
Recent research in graph pebbling has introduced the notion of a
cover pebbling number. Along this same idea, we develop a more
general pebbling function $\pi_{P}^{t} (G)$. This measures the
minimum number of pebbles needed to guarantee that any
distribution of them on $G$ can be transformed via pebbling moves
to a distribution with pebbles on $t$ target vertices.
Furthermore, the $P$ part of the function gives the ability to
change how many pebbles are needed to pebble from one vertex to
another. Bounds on the $\pi$-pebbling function are developed, as
well as its exact value for several families of graphs.
\end{abstract}

\section*{Introduction}

The idea of graph pebbling was first introduced in a paper by
Chung \cite{Start}. If $G$ is a connected graph and $C(G)$ is a
distribution of pebbles onto the vertices of G, a pebbling move
consists of removing two pebbles from a vertex and placing a
pebble on an adjacent vertex. The pebbling number $\pi (G)$ is the
minimum number such that given any configuration $C$ of $\pi (G)$
pebbles, $C$ can be transformed via pebbling moves to a
configuration with a pebble on any target vertex. This question
has been studied extensively and the exact answer is known for a
large set of graphs \cite{Survey}. Recently, several papers have
introduced the concept of a cover pebbling number \cite{Cover}
\cite{Cubes} \cite{Cover2} . This is the minimum number of pebbles
such that any distribution on $G$ can be transformed via pebbling
moves (i.e. pebbled) to end with a distribution where every vertex
has at least one pebble on it.

In this paper, we generalize this concept by asking how many
pebbles are needed to guarantee that a distribution of them can be
pebbled to any $t$ target vertices. This is called the $t$-th
pebbling number. Furthermore, we extend the notion of a pebbling
move by allowing the number of pebbles that are removed at each
vertex to be a function of the vertex.

In this context the first pebbling number is what is classically
called the pebbling number, and the n-th pebbling number is what
is called the cover pebbling number.

In \cite{Cover} the case where $P$ is the standard price function,
$\pi_{P}^{t}$ is called the weighted pebbling number. The paper
defines the $\pi$-function in terms of taking the maximum over a
specific set of weighted cover pebbling numbers, but does not
study it. The problem with using this approach is that previous
research has only calculated positive weighted cover pebbling
numbers, and the maximum must be taken over some non-negative
weighted covering pebbling numbers.

We begin in Section 1 by formalizing the notion of the
$\pi$-pebbling function. Then in Section 2 we prove a theorem
about the cover pebbling number which eliminates degenerate cases
from several proofs. In Section 3 we find appropriate bounds on
the $\pi$-pebbling function for any given graph. These bounds are
then used to calculate the exact value for the complete graph, the
path graph, and star graph in Section 4. In Section 5, we use the
standard price function and study a $\pi$-pebbling sequence.
Section 6 discusses the weighted cover pebbling number and the
cover pebbling theorem. Finally, in Section 7 several open
questions involving the  $\pi$-pebbling function are presented.

\section{Preliminaries}

We begin by formalizing the notion of a configuration of pebbles
on the vertices of a graph G.

\begin{dfn}
A configuration $C(G)$ of $k$ pebbles on the graph G with vertex
set V(G) of size n, is a function $C: V(G) \rightarrow
\mathbb{Z}_{\geq 0}$ such that $\sum_{i=1}^{n} C(v_{i}) = k$. The
set of all configurations on a graph are represented by
$\mathcal{C} (G)$, and the set of all configurations of size $k$
by $\mathcal{C}_{k} (G)$. The value of $C$ at a particular vertex
$i$ is written as $c_{i}$.
\end{dfn}

\begin{dfn}
A price function on a graph G is a function $P:V(G) \rightarrow
\mathbb{Z}_{\geq 2}$. The set of all connected graphs with all
possible price functions is $\mathcal{G}$
\end{dfn}

The value of $P$ at a vertex $i$ is written as either $p_{i}$ or
$p^{i}$. The former is normally used when the vertices are
labelled in some natural way from the graph $G$, and the latter
when we label the vertices such that $p^{i} \leq p^{i+1}$. In
either case, the exact labelling will be made explicit. We now can
define a pebbling move.

\begin{dfn}
A pebbling move is a function $M: C(G) \rightarrow \mathcal{C}(G)$
such that $M(c_{h}) = c_{h} - p_{h}$, $M(c_{k}) = c_{k} + 1$ and
$M(c_{i}) = c_{i}$ for all $i \neq h,k$ where the edge $(v_{h},
v_{k})$ is contained in the edge set $E(G)$. The set
$\mathcal{P}(C)$ is the set of all pebbling moves on the
configuration $C$.
\end{dfn}

One thing to be careful about this definition is that the
resulting M(C) must be in $\mathcal{C}$ and thus can have no
negatively valued vertices. If we want to talk about mapping a
configuration $C$ to another configuration $C'$ via a pebbling
move, we say that $C$ is pebbled to $C'$. With the notion of a
pebbling move, comes the notion of derivability and solvability.
\begin{dfn}
\cite{Condition} We say that the configuration $C'$ is derivable
from $C$ if there are a series of pebbling moves $M_{1}, M_{2},
... M_{n}$ such that $M_{n} \circ ... M_{2} \circ M_{1} (C) = C'$.
\end{dfn}

\begin{dfn} A configuration $C$ is said to cover a subset of $V(G)$
if $c_{i}$ is non-zero on the entire subset. A configuration $C$
is t-solvable \footnote[1]{This is not to be confused with the
definition of t-solvability given in \cite{Survey}}, if given a
set of t vertices, there exists a configuration that covers those
t vertices and is derivable from $C$.
\end{dfn}

Now we finally can define the notion of a $\pi$-pebbling function.

\begin{dfn}
The $\pi$-pebbling function is a map $\pi: \mathcal{G} \times
[1,n] \cap \mathbb{Z} \rightarrow \mathbb{Z}_{\geq 1}$. The value
of $\pi (G, P, t)$, denoted by $\pi_{P}^{t} (G)$, is the minimum
number $k$ such that any configuration C contained in
$\mathcal{C}_{k} (G)$ with price function $p$ is t-solvable.
\end{dfn}

\section{Theorem on the Cover Pebbling Number}

One difficulty in dealing with the general pebbling function, is
that oftentimes the $n$-th pebbling number, the cover pebbling
number, does not fit directly into the formula for the other
pebbling numbers. Instead of dealing with this case by case, there
is a nice theorem which relates $\pi_{P}^{n-1}$ and $\pi_{P}^{n}$.

\begin{thm}
For any graph $G$ with $n$ vertices, $\pi_{P}^{n-1} (G) + 1 =
\pi_{P}^{n} (G)$
\end{thm}

Proof: Let $k = \pi_{P}^{n-1} (G)$, and let $C$ be a configuration
of $k+1$ pebbles on $G$. We show that there is a way to pebble $C$
to a covering of $G$. Pick a vertex $v_{1}$ that has at least one
pebble from $C$, and mark one of the pebbles on this vertex. Now,
since $k = \pi_{P}^{n-1} (G)$, there is a way to pebble $C$ onto
$G - v_{1}$ without moving the marked pebble. This results in a
covering of $G$. The fact that $\pi_{P}^{n} (G) > \pi_{P}^{n-1}
(G)$ is obvious. $\mathbb{Q.E.D.}$

\medskip

When the cover pebbling number does not fit into the general
pattern, it will be noted. This theorem can then be used to get it
from the given formula.

\section{Bounds on the $\pi$-Pebbling Function}

We begin by finding the obvious general upper and lower bounds on
the $\pi$-pebbling function. The first lower bound comes from
assuming that the pebbles are spread out over all but t of the
vertices.

\begin{thm}
Let $d$ be the diameter of $G$ and the vertices of $G$ be numbered
such that $p^{i} \leq p^{i+1}$, then for $n \neq t$
\begin{displaymath}
\sum_{i=t+1}^{n} (p^{i}-1)+p^{n} (t-1) +1 \leq \pi_{C}^{t} (G)
\end{displaymath}
\end{thm}

Proof: Assume $(p^{i}-1)$ pebbles are placed on their respective
vertices for all $i$ such that $t+1 \leq i \leq n-1$, and
$p^{n}t-1$ pebbles are placed on the $p^{n}$th vertex. Then there
are $t$ unpebbled vertices and $\sum_{i=t+1}^{n-1} (p^{i}-1)+p^{n}
t - 1$ pebbles. The only vertex that can be pebbled from is the
$p^{n}$th, and there are only enough pebbles to pebble to $t-1$
vertices, so we have $\pi_{P}^{t} (G) \geq \sum_{i=t+1}^{n-1}
(p^{i}-1)+p^{n} t$. $\mathbb{Q.E.D.}$

\medskip

The second lower bound comes from assuming all the pebbles are
placed on the same vertex. This requires two different results,
one for $t \leq d$ and one for $t > d$.

\begin{thm}
Let the diameter $d$ of $G$ be greater than or equal to t, and
$\Gamma$ be a path of length $d$ in $G$. Furthermore, let the path
$\Gamma$ be numbered such that $p_{i}$ is adjacent to $p_{i+1}$
for all $i$ and under this restriction the product
$p_{1}p_{2}...p_{n-d}$ is maximal. Then for $t < d$,
\begin{displaymath}
\sum_{j=1}^{t} (\prod_{i=1}^{n-j} p_{j}) \leq \pi_{P}^{t} (G)
\end{displaymath}
Now if we have $t$ is greater than $d$ and $n \neq t$ then,
\begin{displaymath}
\sum_{j=1}^{t} (\prod_{i=1}^{n-j} p_{j}) + (p_{1} p_{2}) (d-t)
\leq \pi_{P}^{t} (G)
\end{displaymath}
\end{thm}

Proof: $(t \leq d)$ Since the graph has diameter $d$, this means
that there is a path in $G$ with length $d$. Consider
$\sum_{j=1}^{t} (\prod_{i=2}^{n-j+1} p^{j})-1$ pebbles placed on
$p_{1}$. Even if the price function on the graph is as low as
possible, there are not enough pebbles to pebble to the farthest
$t$ vertices on the path.

$(t > d)$ Using the same argument, place $\sum_{j=1}^{d}
(\prod_{i=2}^{n-j+1} p^{j}) + (p^{1} p^{2}) (d-t)-1$ pebbles on an
endpoint of a length $d$ path. At the least it will take
$\sum_{j=1}^{d} (\prod_{i=2}^{n-j+1} p^{j})$ pebbles to fill up
the other $d$ vertices in the path. Then, in order to fill up
$t-d$ vertices not on the path, the least number of pebbles needed
would be $(p^{1} p^{2}) (d-t)$ since none of these vertices could
be directly adjacent to the initial vertex or else there would be
a path of length $d+1$. $\mathbb{Q.E.D}$

\medskip

There is also a corresponding higher bound on the $\pi$-pebbling
function. This comes from putting all the pebbles on one vertex
and assuming the price function is as high as possible from this
vertex.

\begin{thm}
Let $d$ be the diameter of a graph $G$, then
\begin{displaymath}
\pi_{P}^{t} (G) \leq t [(\prod_{i=n-(d-1)}^{n} (p^{i}) -
1)(n-1)+1]
\end{displaymath}
\end{thm}

Proof: Looking at the case of $t=1$, assume $
[(\prod_{i=n-(d-1)}^{n} (p^{i}) - 1)(n-1)+1]$ pebbles are
distributed on the graph $G$. Then either every vertex has a
pebble on it or some vertex has $\prod_{i=n-(d-1)}^{n} (p^{i})$
pebbles on it by the pigeonhole principle. This implies the
theorem for $t=1$. Now notice that since $\pi_{P}^{t} (G)$ is the
number of pebbles needed to pebble to any vertex, we have the
inequality $\pi_{P}^{t}(G) \leq t \pi_{p}^{1} (G)$. This finishes
the proof. $\mathbb{Q.E.D.}$

\medskip

These three inequalities are nice because in the case of $t=1$ and
$p_{i} = 2$ for all $i$ they all collapse to the inequalities for
the classical pebbling number given in \cite{Survey}.

\begin{cor}
Let $G$ be a graph with diameter $d$, then $max[n,2^{d}] \leq
\pi^{1} (G) \leq (2^{d} -1)(n-1) +1$.
\end{cor}

This inequality immediately gives the first pebbling number for
$K_{n}$ and $P_{n}$. These two graphs in turn show the sharpness
for the case of $t=1$ of all three bounds. It is natural to ask
whether or not the general bounds are sharp for all values of $t$.
As will be shown in the next section, the complete graph and path
graph make Theorem 8 and Theorem 9 sharp respectively for all
values of $t$ and $n$. Using Theorem 7, it is not hard however to
show that the inequality it Theorem 10 is equal only when $n$ is
1.

\section{Complete Graphs, Path Graphs, and Star Graphs}

The first family of graphs we look at are the complete graphs.
These are not too complicated since any set of $t$ vertices are
indistinguishable from another set of $t$ vertices. The complete
graphs are also nice because they will show that the lower bound
given in Theorem 8 is sharp.

\begin{thm}
If the vertices of $K_{n}$ are numbered such that $p^{i} \leq
p^{i+1}$, then for $n \neq t$
\begin{displaymath}
\pi_{P}^{t} (K_{n}) = \sum_{i=t+1}^{n} (p^{i}-1)+p^{n} (t-1) +1.
\end{displaymath}
\end{thm}

Proof: Let $C$ be a configuration of $\sum_{i=t+1}^{n}
(p^{i}-1)+p^{n} (t-1) +1$ pebbles on $K_{n}$. Assume that $t$
target vertices have been selected. We can assume that there is at
least one target vertex $v_{1}$ such that $c_{1}=0$, and then
there are $\sum_{i=t+1}^{n-1} (p^{i}-1)+p^{n} t$ pebbles on $n-1$
vertices. The pigeonhole principle says that there exists a vertex
such that $c_{i} \geq p_{i}$. This implies that there is a way to
get a pebble to $v_{1}$. Now if we continue to fill up empty
target vertices in this way, the pigeonhole principle says that
this can continue for at least $t$ steps. The only question is
whether or it is possible that in the process of coving a target
vertex, another target vertex becomes uncovered. It turns out that
this can only happen if the target vertex was covered by the
original configuration, and not by pebbling. This is because in
each pebbling step the chosen target vertex ends up with only 1
pebble, so there is no way to pebble off of it in future pebbling
moves. Therefore, using this algorithm until there are no
uncovered target vertices will take at most $t$ steps, and thus
there are enough pebbles to carry it out. Theorem 8 gives the
lower bound and completes the result. $\mathbb{Q.E.D.}$

\medskip

While the path graph is not symmetric, it turns out that using a
convenient numbering of the vertices the pebbling number is rather
simple. The path graph also shows that Theorem 9 is sharp, which
is natural since it involves putting all the pebbles on the end of
a maximal path in $G$.

\begin{thm}
If the vertices on the path of n elements are labelled such that
$v_{i}$ is adjacent to $v_{i+1}$ for all i and under this
restriction $p_{1} p_{2} ... p_{t}$ is as large as possible, then
for $n \neq t$
\begin{displaymath}
\pi_{P}^{t} (P_{n}) = \sum_{j=1}^{t} (\prod_{i=1}^{n-j} p_{i})
\end{displaymath}
\end{thm}

Proof: Let D be a distribution of $\sum_{j=1}^{t}
(\prod_{i=1}^{n-j} P_{i})$ pebbles on $P_{n}$. If all the vertices
have pebbles, then we are done, so assume there exists a vertex
without any pebbles. Then there are $\sum_{j=1}^{t}
(\prod_{i=1}^{n-j} p_{i})$ pebbles on $n-1$ vertices. Using the
pigeonhole principle, it can be shown that there is a way to
pebble to any open vertex. Now assume we have a distribution of
size $\sum_{j=1}^{t} (\prod_{i=1}^{n-j} p_{i})$ with at least $k
\leq t$ vertices without pebbles and we can pebble to any $k-1$ of
them. If we pebble to all but the lowest numbered vertex, then
there are at the least $\sum_{j=k}^{t} (\prod_{i=1}^{n-j} p_{i})$
pebbles left on the other $n-(k-1)$ vertices and the pigeonhole
principle again shows that these can be pebbled to the $k$-th
vertex. This upper bound is the same as the lower bound in Theorem
9, and thus the result follows. $\mathbb{Q.E.D.}$

\medskip

Now we are ready to look at the star graph $S_{n}$. The star graph
consists of a central node, connected to $n$ degree one vertices.
This graph is not regular and has to be dealt with in two cases.

\begin{thm}
If the vertices of the star graph are labelled such that $p_{0}$
is the price of the center vertex and $p^{i} \leq p^{i+1}$ for all
other $i$, then for $(t < n)$
\begin{displaymath}
\pi_{P}^{t} (S_{n}) = \sum_{i=t+1}^{n-1} (p^{i} - 1) + t p_{0}
p_{n}
\end{displaymath}
for $(t = n)$,
\begin{displaymath}
 \pi_{P}^{t} (S_{n}) = n p_{0} p^{n} + p^{n}
\end{displaymath}
\end{thm}

Proof:($t < n$) Assume $\sum_{i=t+1}^{n-1} (p^{i} - 1) + t
p_{0}p_{n} -1$ pebbles are placed on the vertex $p^{n}$ and
$(p^{i}-1)$ are put on there respective vertices for for $t+1 \leq
i \leq n-1$. There are then $t$ empty outer vertices. The only
vertex that can be pebbled from is $p_{n}$, and there are not
enough pebbles to pebble to all $t$ empty vertices. Therefore
$\sum_{i=t+1}^{n-1} (p^{i} - 1) + t p_{0}p_{n} \leq \pi_{P}^{t}
(S_{n})$.

Now assume that there is a configuration of $\sum_{i=t+1}^{n-1}
(p^{i} - 1) + t p_{0}p_{n}$ pebbles on $S_{n}$. We assume that
there are no pebbles on the central vertex. Assume there are $k
\leq t$ empty outer vertices. Then there are $\sum_{i=t+1}^{n-1}
(p^{i} - 1) + t p_{0}p_{n}$ on $n-k$ vertices. The pigeonhole
principle says that as long as $k \leq t$ there is a way to pebble
$p_{0} t$ pebbles on to the central node, and thus a way to cover
the $k$ empty outer vertices. The case of $c_{0} \neq 0$ or the
central node being a target vertex is obvious given the above
proof.

($t = n$) Assume $n p_{0} p^{n} + p^{n} -1$ pebbles are placed on
the vertex $p^{n}$. There are then $n$ empty outer vertices. The
only vertex that can be pebbled from is $p_{n}$, and there are not
enough pebbles to pebble to all $t$ empty vertices. Therefore $n
p_{0} p^{n} + p^{n} \leq \pi_{P}^{t} (S_{n})$. Showing the other
inequality is the same as above. $\mathbb{Q.E.D.}$

\section{$\pi$-Pebbling Ratio Series}

In \cite{Cover} a covering ratio of a graph $G$ \footnote[2]{In
this section we assume that G has the standard price function
$p_{i} = 2$ for all $i$ and omit the P from the $\pi$-pebbling
function} is defined as $\pi^{n} (G) / \pi^{1} (G)$. Since we have
a whole set of pebbling numbers, the natural generalization of
this is a sequence.

\begin{dfn}
The $\pi$-pebbling ratio series of the graph $G$ is the series $(
\pi^{i+1} (G) / \pi^{i} (G))_{i=1}^{n-1}$. We denote it by $\rho$
and the $i$-th term by $\rho_{i}$.
\end{dfn}

If we have an entire family of graphs, such as the complete
graphs, there is a natural way to define an infinite sequence on
the whole family.

\begin{dfn}
Let $\mathcal{F}$ be a family of graphs such that for every
natural number n there is an unique graph in the family with $n$
vertices. Then the pebbling ratio sequence of the family is the
sequence $(\rho_{t} (F_{i}))_{i=t}^{\infty}$ where $F_{n}$ is the
unique representative with n vertices. We denote this sequence as
$\alpha^{t} (\mathcal{F})$.
\end{dfn}

The nice thing about the pebbling ratio sequence of a family of
graphs is that it is infinite, and so we can talk about when it
converges and what it converges to.

Now using our results on complete graphs, path graphs, and star
graphs, we can determine what the ratio sequence of these families
converge too for various values of t.

\begin{thm}
Let $\mathcal{P}$ be the family of paths on n vertices. Then
$\alpha^{t} (\mathcal{P})$ is a constant series and the constant
is $\frac{2^{t+2} - 1}{2(2^{t} -1)}$.
\end{thm}

Proof: We see from Theorem 13, that $\pi^{t} (P_{n})= 2^{n} -
2^{n-t}$. Therefore
\begin{displaymath}
\alpha^{t}_{n} (\mathcal{P}) = \frac{ 2^{n} - 2^{n-(t+1)}}{ 2^{n}
- 2^{n-t}}
\end{displaymath}
By multiplying the top and bottom by 2 we get
\begin{displaymath}
\alpha^{t}_{n} (\mathcal{P}) = \frac{ 2^{n+1} - 2^{n+1-(t+1)}}{
2^{n+1} - 2^{n+1-t}} = \alpha^{t}_{n+1} (\mathcal{P})
\end{displaymath}
And therefore the series is constant. To get the constant set n=1
and simplify since the series is constant. $\mathbb{Q.E.D.}$

\begin{thm}
Let $\mathcal{K}$ be the family of complete graphs on n vertices.
Then $\alpha^{t} (\mathcal{K})$ converges to 2 for all values of
t.
\end{thm}

Proof: We see from Theorem 12, that $\pi^{t} (K_{n})=
n-t+2^{n}(t-1) +1$. Therefore we can take the limit
\begin{displaymath}
\lim_{n \rightarrow \infty} \frac{n+1-t+2^{n+1}(t-1)
+1}{n-t+2^{n}(t-1) +1}
\end{displaymath}
Taking the limit yields the result of 2. $\mathbb{Q.E.D.}$

\begin{thm}
Let $\mathcal{S}$ be the family of star graphs on n vertices. Then
$\alpha^{t} (\mathcal{S})$ converges to 1 for all values of t.
\end{thm}

Proof: We see from Theorem 14, that $\pi^{t} (S_{n})= 3t+n-1$.
Therefore we can take the limit
\begin{displaymath}
\lim_{n \rightarrow \infty} \frac{3t+n}{3t+n-1}
\end{displaymath}
Dividing by n and taking the limit yields the result.
$\mathbb{Q.E.D.}$

Now that we have shown that these sequences converge, we define
the secondary pebbling ratio sequence.

\begin{dfn}
If $\alpha^{t} (\mathcal{F})$ converges for all t for some family
of graphs, we can define the secondary (or Beta) pebbling ratio
sequence as the sequence $\beta (\mathcal{F}) = (\alpha^{t}
(\mathcal{F}))_{t=1}^{\infty}$.
\end{dfn}

The $\beta$ pebbling sequence can be thought of as the $\rho$
pebbling sequence of $F_{\infty}$. Now we can calculate $\beta
(\mathcal{F})$ for these three families and get $\beta
(\mathcal{P})$ converges to 2, $\beta (\mathcal{K})$ converges to
2, and $\beta (\mathcal{S})$ converges to 1. What this means is
the given a large enough n and a large enough t, in order to be
able to pebble to $t+1$ vertices takes about twice as many pebbles
in the complete graph or the path graph but takes about the same
number in the star graph.

These two sequences appear to have an extraordinary amount of
structure. In order to fully understand them more $\pi$-pebbling
numbers must be constructed for more families of graphs. One
approach to doing this is described in the next section. One thing
to note is that some families of graphs, such as the family of all
trees or all odd cycles, do not have the property that there is
one unique graph on n vertices for every natural number n. In the
case of odd cycles when existence is the problem, one could simply
define the pebbling ratio sequence in the normal way, but simply
skip terms that are not defined. In the case where uniqueness is a
problem, such as trees, the fix is more difficult.

\section{Weighted Cover Pebbling Number}

A natural way to generalize our results on the $\pi$-pebbling
function would be to ask what is the minimum number of pebbles
needed to guarantee that a sequence of pebbling moves can pebble
to a configuration with a minimum number of pebbles on each
vertex. This idea is worked out in \cite{Cover}.

\begin{dfn}
A weight function is a map: $W: V(G) \rightarrow \mathbb{Z}_{\geq
0}$. Its value at a vertex $i$ is denoted by $w_{i}$. A weight
function is said to be positive if $w_{i} \geq 1$ for all $i$. The
size of W is the sum of its values over all vertices
\end{dfn}

\begin{dfn}
The weighted cover pebbling number $\gamma_{W} (G)$ is the minimum
number k such that any $C \in \mathcal{C}_k$ can be pebbled to a
configuration $C'$ such that $c'_{i} \geq w_{i}$.
\end{dfn}

We see that by taking the max over all weighted cover pebbling
numbers of size $k$ on a graph $G$, we would get $\pi^{t} (G)$.
But while the weighted cover pebbling number has been studied for
the case of a positive weight, it has not been calculated for
non-positive weight functions. This is because all positive weight
functions share a nice property that makes their cover pebble
number easy to calculate.

\begin{thm} \cite{Theorem} (Cover Pebbling Theorem) If $\gamma_{W} (G) = k$, then there
exists a simple configuration of size $k-1$ that is not solvable
for the weight function $W$.
\end{thm}

We have seen from the complete graph and the star graph that this
is not the case for non-positive weight functions. A very
important result in studying the $\pi$-pebbling function would be
to find an appropriate generalization of the Cover Pebbling
Theorem. This would allow one to calculate the number for more
complicated graphs like wheels without too much difficulty. We
provide a conjecture on what a possible generalization might be.

\begin{con}
Let $\pi^{t} (G) = k$. Then there exists a non-t-solvable
configuration $C$ of $k-1$ pebbles on G such that $c_{i}$ is
either $0$, or $p_{i}-1$ except for possibly one vertex $c_{r}$.
For the standard price function, the point $r$ should be an
element that has another point $s$ a distance $d$ away, where $d$
is the diameter of the graph, and there should exist a minimal
path from $r$ to one such $s$ that has no pebbles on it other than
those on $r$.
\end{con}

\section{Open Questions}

There are so many unknown things about the $\pi$-pebbling function
that we present only a couple problems that the author sees to be
the most important.

\begin{opn}
Prove some variation on Corollary 24.
\end{opn}

\begin{opn}
Determine the weighted pebbling number for weights that are not
necessarily positive.
\end{opn}

\begin{opn}
Determine when the primary and secondary pebbling ratio series
converge, and determine what values they can converge to (natural
numbers?).
\end{opn}

\begin{opn}
Using the full $\pi$-pebbling function, determine bounds on the
pebbling threshold of various families of graphs.
\end{opn}

\section*{Acknowledgements}

The author wishes to thank Steve Fisk for introducing him to the
area of graph pebbling.



\begin{thebibliography}{77} 

\small

\bibitem{Start} F. R. K. Chung {\it Pebbling in Hypercubes,} SIAM J. Disc. Math (1989), 467--472.

\bibitem{Cover} Betsy Crull, Tammy Cundiff, Paul Feltman, Glenn H.
Hurlbert, Lara Pudwell, Zsuzsuanna Szaniszlo, Zsolt Tuza {\it The
Cover Pebbling Number of Graphs,} (2005), Math ArXiv
math.CO/0409368.

\bibitem{Survey} Glenn H.Hurlbert,
``A Survey of Graph Pebbling,'' {\it Congressus Numerantium,} {\bf
139} (199): 41--64.

\bibitem{Cubes} Glenn H. Hurlbert, Benjamin Munyan,
{\it Cover Pebbling Hypercubes,} Math ArXiv math.CO/0409321.

\bibitem{Theorem} Jonas Sjostrand, The Cover Pebbling Theorem,
Math ArXiv, math.co/0410129

\bibitem{Condition} Annal Vuong, M. Ian Wyckoff, {\it Conditions of Weighted Cover Pebblings of Graphs,} Math ArXiv math.CO/0410410v1.

\bibitem{Cover2} Nathaniel G. Watson, Carl. R. Yerger
{\it The Cover Pebbling Numbers and Bounds for Certain Families of
Graphs,} (New York: Addison, 1994).
\end{thebibliography}
\end{document}